\documentclass[12pt,a4paper]{article}
\usepackage{float}
\usepackage{amsmath}
\usepackage{amsfonts}
\usepackage{amssymb}
\usepackage{graphicx}%
\renewcommand{\baselinestretch}{1}
\newtheorem{theorem}{Theorem}

\newtheorem{lemma}{Lemma}

\newtheorem{proposition}{Proposition}
\newtheorem{remark}{Remark}

\newenvironment{proof}[1][Proof]{\noindent\textbf{#1.} }{\ \rule{0.5em}{0.5em}}

\let \a = \alpha
\let \b = \beta
\let \l = \lambda
\let \g = \gamma

\def \R {{\mathbb R}}

\begin{document}

\title{\bf On a general similarity boundary layer equation}
\author{BERNARD BRIGHI$^a$  and JEAN-DAVID HOERNEL$^{b,*}$ }
\date{}
\maketitle

\begin{center}
Universit\'e de Haute-Alsace 
\vskip 0,1cm
Laboratoire de Math\'ematiques, Informatique et Applications
\vskip 0,1cm
4 rue des fr\`eres Lumi\`ere, 68093 MULHOUSE (France)
\end{center}

\begin{abstract}
In this paper we are concerned with the solutions of the differential equation $f^{\prime\prime\prime}+ff^{\prime\prime}+g(f^{\prime})=0$ on $[0,\infty)$, satisfying the boundary conditions $f(0)=\alpha$, $f'(0)=\beta\geq 0$, $f'(\infty)=\l$, and where $g$ is some given continuous function.
This general boundary value problem includes the Falkner-Skan case, and can be applied, for example, to free or mixed convection in porous medium, or flow adjacent to stretching walls in the context of boundary layer approximation.
Under some assumptions on the function $g$, we prove existence and uniqueness of a concave or a convex solution. We also give some results about nonexistence and asymptotic behaviour of the solution.
\end{abstract}



\renewcommand{\baselinestretch}{1}
\footnotetext{\hspace{-0.8cm}$^*$Corresponding author.}
\footnotetext{\hspace{-0.8cm}E-mail addresses: bernard.brighi@uha.fr, j-d.hoernel@wanadoo.fr}

\section{Introduction}
We consider the following third order non-linear autonomous differential
equation
\begin{equation}
f^{\prime\prime\prime}+ff^{\prime\prime}+g(f^{\prime})=0 \label{eq}
\end{equation}
with the boundary conditions
\begin{eqnarray}
f(0)=\alpha, \label{c1} \\
f'(0)=\beta, \label{c2} \\
f'(\infty)=\l \label{c3}
\end{eqnarray}
where $\a\in \R$, $\b\in \R_+$, $\l \in \R$ and $f'(\infty):=\underset{t\rightarrow\infty}{\lim} f'(t)$. We also assume that the given function $g$ is locally Lipschitz on some interval $J$ containing $\b$ and $\l$.

\newpage


In the litterature, problem (\ref{eq})-(\ref{c3}) with suitable $g$ and $\l$ arises in many fields of application such as free or mixed convection in a fluid saturated porous medium near a semi or double infinite wall in the framework of boundary layer approximation, high frequency excitation of liquid metal, stretching walls,...  

The problems of free convection, stretching walls and high frequency excitation of liquid metal corresponds to the function $g$ given by $g(x)=\frac{2m}{m+1}(-x)x$ and to $\l=0$. There are many papers in connection with those such as \cite{banks1}, \cite{banks},  \cite{brighi04},  \cite{pop1}, \cite{pop}, \cite{cheng}, \cite{crane}, \cite{ene}, \cite{gup}, \cite{ing}, \cite{mag}, \cite{mag1}, \cite{merk}, \cite{mof} for the physical point of view or numerical computations,  \cite{brighicr}, \cite{brighi02}, \cite{brighi01},  \cite{heat_flux}, \cite{equiv}, \cite{brighisari}, \cite{guedda}, \cite{guedda1} for the mathematical analysis and \cite{gaeta} for a survey.

The Falkner-Skan equation, arising in the study of two dimensional flow of a slightly viscous incompressible fluid past a wedge of angle $\pi m$ under the assumptions of boundary-layer theory, is obtained for $g(x)=m(1-x)(1+x)$ and $\l=1$. This famous equation has been widely studied, see for example \cite{coppel}, \cite{falk}, \cite{hart}, \cite{hast}, \cite{hast1} and the references therein for a survey, and \cite{ish1}, \cite{wang}, \cite{yang}, \cite{yang2} for more recent investigations.

The mixed convection case corresponds to $g(x)=\frac{2m}{m+1}(1-x)x$ and  $\l=1$. This problems appears recently in \cite{aly} and \cite{nazar}. Some first theoretical results about this equation can be found in \cite{aml} and \cite{guedda2}.

The Blasius problem, corresponding to $g=0$, is a particular case of all the previous situations and the first historic case in which an equation of the form (\ref{eq}) appears.
 This well known problem, that arises in \cite{bla} at the begining of the previous century, has been studied in a lot of papers. For more details, we refer to \cite{brighi03}, \cite{coppel} and \cite{hart} and the references therein.


Finally, let us notice that a first generalization of some of the previous equations can be found in \cite{utz}. The author considers problem (\ref{eq})-(\ref{c3}) with $\a\geq 0$, $\l\geq \b>0$ and functions $g$ such that $g(x)=\hat g(x^2)$ where $\hat g$ is assumed to be positive and monotone decreasing on $[\b,\l)$ and $\hat g(\l^2)=0$. Under these hypotheses, he proves that there exists one and only one convex solution of this problem.

\begin{remark}
Let $a>0$, and consider the differential equation $u'''+auu''+h(u')=0$ where $h$ is some given function $h$. By setting $u(t)=\frac{1}{\sqrt{a}}f(\sqrt{a}t)$ then the equation in $u$ reduces to $f'''+ff''+g(f')=0$ with $g(x)=\frac{1}{a}h(x)$. 


\end{remark}

\section{Preliminary results}

First of all, let us remark that if $f$ satisfies equation (\ref{eq}) on some interval $I$ and if we denote by $F$ any anti-derivative of $f$ on $I$, then we have
\begin{equation}
\left (f''e^F\right )'=-g(f')e^F.\label{exp}
\end{equation}
This, in particular implies that the concavity of $f$ is related to the changes of sign of $g$.

We will need some Lemmas concerning the solutions of (\ref{eq}).
\begin{lemma}\label{vanish}
If $g(\mu)=0$ and if $f$ is a solution of {\rm (\ref{eq})} on some interval $I$ such that there exists a point $t_0\in I$ verifying $f''(t_0)=0$ and $f'(t_0)=\mu$, then $f''(t)=0$ for every $t\in I$.
\end{lemma}
\begin{proof}
Let $f$ be a solution of (\ref{eq}) on $I$ such that $f''(t_0)=0$ and Ê$f'(t_0)=\mu$ for some $t_0\in I$. Since the function $r(t)=\mu(t-t_0)+f(t_0)$ is a solution of (\ref{eq}) such that $r(t_0)=f(t_0)$, $r'(t_0)=f'(t_0)$ and $r''(t_0)=f''(t_0)$, we get $r=f$ and $f''\equiv 0$.
\end{proof}

\begin{lemma}\label{f2} Let $f$ be a solution of $(\ref{eq})$ on some interval
$[t_0,\infty)$, such that $f'(t)\to l \in\R$ as $t\to\infty$. If moreover $f$ is of constant sign at infinity, then, we have
$$\lim_{t\to\infty}f''(t)=0.$$
\end{lemma}
\begin{proof} First of all, let us remark that since $f'(t)$ has a finite limit as
$t\to\infty$, then 
\begin{equation}
\liminf_{t\to\infty}f''(t)^2=0.\label{i}
\end{equation}
Multiplying $(\ref{eq})$ by $f''$ and integrating on $[t_0,t]$, we get 
\begin{equation}
\frac{1}{2}f''(t)^2-\frac{1}{2}f''(t_0)^2+\int_{t_0}^tf(s)f''(s)^2ds
+G(f'(t))-G(f'(t_0))=0,\label{ii}
\end{equation} 
where we denoted by $G$ any anti-derivative of $g$. 
As $f$ is of constant sign at infinity, it follows that the integral in $(\ref{ii})$, and
thus $f''(t)^2$ too, have limits as $t\to\infty$. From $(\ref{i})$ we get the result.
\end{proof}

\begin{remark} If $l \not=0$ we have $f(t)\sim l t$ as $t\to\infty$ and $f$ is of constant sign at infinity.
If $l=0$ and if $f$ is either concave or convex at infinity, then again $f$ is of constant sign
at infinity. From $(\ref{exp})$ it is the case, for example, if $g$ is of constant sign in a neighbourhood of $0$.
\end{remark}

The following  Lemma shows that to expect a solution of $(\ref{eq})$-$(\ref{c3})$, we must assume that the function $g$ vanishes at the point $\l$. 

\begin{lemma} \label{lim} Let $f$ be a solution of $(\ref{eq})$ on some interval $[t_0,\infty)$, such that $f'(t)\to l\in \R$ as $t\to\infty$. Then $g(l)=0$.
\end{lemma}
\begin{proof}
Let us suppose that $2c=-g(l)>0$. There exists $t_1>t_0$ such that $-g(f'(t))>c$ for $t>t_1$ and from (\ref{exp}) we have
$$\left ( f''e^{F}\right )'>ce^{F}$$
on $[t_1,\infty)$. This means that $f''$ cannot vanish more than once and thus $f$ is concave or convex at infinity.
\begin{itemize}
\item Assume now that $f$ is bounded. Then, it follows from $(\ref{eq})$ and Lemma \ref{f2} that $f'''(t)\to 2c$ as $t\to\infty$ and we have a contradiction.
\item Assume next that $f$ is unbounded. Then $\vert f(t)\vert\to\infty$ as $t\to\infty$ and thus there exists $t_2>t_1$ such that $f'/f>0$ on $[t_2,\infty)$. From $(\ref{eq})$ we get 
\begin{equation}
\frac{f'''f'}{f}+f''f'>c\frac{f'}{f} \quad \text{on} \quad [t_2,\infty).\label{iii}
\end{equation}
Integrating we obtain
$$\int_{t_2}^t\frac{f'''(s)f'(s)}{f(s)}ds+\frac{1}{2}f'(t)^2-\frac{1}{2}f'(t_2)^2\geq c\left(\ln\vert f(t)\vert-\ln\vert f(t_2)\vert\right).$$
It follows that 
\begin{equation}
\int_{t_2}^\infty\frac{f'''(s)f'(s)}{f(s)}ds=\infty.\label{iv}
\end{equation}
But, we have
\begin{align*}
\int_{t_2}^t\frac{f'''(s)f'(s)}{f(s)}ds= &\ \frac{f''(t)f'(t)}{f(t)}-\frac{f''(t_2)f'(t_2)}{f(t_2)} \cr
&\ -\int_{t_2}^t \frac{f''(s)^2}{f(s)}ds+\int_{t_2}^t \frac{f''(s)f'(s)^2}{f(s)^2}ds
\end{align*}
which leads to a contradiction with $(\ref{iv})$, since the integrals in the right hand side have
finite limits as $t\to\infty$. 

For $c<0$, same arguments give also a contradiction. The proof is now complete.
\end{itemize}
\vspace{-3mm}
\end{proof}

\begin{remark}
If $l \not = 0$ we can have a much simpler proof. Indeed, in this case we have $f'(t)\sim l$ and $f(t) \sim lt$ as $t \to \infty$. Integrating $(\ref{eq})$ on $[t_0,t]$ we get
\begin{align*}
f''(t)-f''(t_0)-f(t_0)f'(t_0)&=-f(t)f'(t)+\int_{t_0}^t f'(s)^2ds-\int_{t_0}^tg(f'(s))ds \cr
&=-l^2t(1+o(1))+l^2t(1+o(1))-g(l)t(1+o(1)) \cr
&=-g(l)t+o(t)
\end{align*}
which is a contradiction since $f''(t)\to 0$ as $t\to \infty$ by Lemma $\ref{f2}$.
\end{remark}

\begin{remark} Solution, for which first derivative does not have a finite limit, does exist. 
For example
\begin{itemize} \item For any $a\in\R^*$ and any $b\in\R$ the function $f$ defined by
$$f(t)=at^2+bt+\frac{b^2-1}{4a}$$
is a solution of $(\ref{eq})$ with $g(x)=\frac{1}{2}(1-x^2)$ and $f'(t)\to \pm\infty$ as $t\to \infty$ $($Falkner-Skan with $m=\frac{1}{2})$.
\item If $g(x)=-x^2+x+1$, then $f(t)=\sin t$ is a solution of $(\ref{eq})$ for which $f'$ does not have a limit at infinity.
\end{itemize}
\end{remark}


In order to get solutions of (\ref{eq})-(\ref{c3}) for given $\alpha\in\mathbb R$, $\beta \in\mathbb R_+$ and $\l\in\mathbb R$, we will consider the  initial value problem 
\begin{equation}
\left \{ \begin{array}
[c]{lll}
f'''+ff''+g(f')=0,\\ 
f(0) = \a,\\
f'(0) = \b,\\
f''(0) = \g
\end{array} \right . \label{initial}
\end{equation}
and use a shooting technique on the parameter $\gamma$. We will denote by $f_\gamma$ its solution and by $[0,T_\gamma)$ its right maximal interval of existence. 
Integrating (\ref{eq}) on $[0,t]$ for $0<t<T_\gamma$, we obtain the useful identity
\begin{equation}
f''_\gamma(t)=\g -f_\g(t)f'_\g(t)+\a\b+\int_0^t \left (f_\g'(s)^2-g(f_\g'(s))\right )ds. \label{i1}
\end{equation}

\begin{remark}
Looking at $(\ref{i1})$ we see that if we take $g(x)=x^2$ the integral in the right-hand side vanishes. Then, integrating on $[0,t]$ we obtain
$$f'(t)+\frac{1}{2}f(t)^2=(\g+\a\b)t+\frac{\a^2}{2}+\b$$
and choosing $\g=-\a\b$ we have that for $\b>-\a^2/2$ the function
$$f(t)=\frac{1}{-\frac{1}{2d}+\frac{\a+d}{2\a d-2d^2}e^{dt}}+d$$
with $d=\sqrt{\a^2+2\b}$ is a bounded solution of the problem {\rm (\ref{eq})-(\ref{c3})} for $\l=0$.
\end{remark}

\begin{remark}
Let us take a look at the case $\beta=\l>0$. If $g(\l)=0$, then the function $f_0(t)=\l t+\alpha$ is  a solution of $(\ref{eq})$-$(\ref{c3})$. Without additional hypotheses on $g$, we cannot say anythings about uniqueness. However, if we assume, for example, that $g(x)<0$ for $x>\l$ and $g(x)>0$ for $x<\l$, then $f_0$ is the unique solution of $(\ref{eq})$-$(\ref{c3})$. Indeed, let $f_\g$ be another solution of $(\ref{eq})$-$(\ref{c3})$ such that $\gamma>0$. Since $f_\g'(0)=f_\g'(\infty)=\l$, there exists $t_0>0$ such that $f_\g'(t_0)>\l$, $f_\g''(t_0)=0$ and $f_\g'''(t_0)\leq 0$. From $(\ref{eq})$ we obtain $f_\g'''(t_0)=-g(f_\g'(t_0))>0$ and thus a contradiction. 
If $\gamma<0$, the same approach leads again to a contradiction.
\end{remark}

In the following, we will focus first on the concave solutions and next on the convex solutions of (\ref{eq})-(\ref{c3}) with functions $g$ such that $g(\l)=0$. As seen in Lemma \ref{lim}, this hypothesis  is necessary to realize the condition $f'(t)\to \lambda$ as $t\to \infty$. In addition, we will assume that some condition on the sign of $g$ is satisfied between $\b$ and $\l$ in order to get existence and uniqueness of a concave solution (when $\l<\b$) or a convex solution (when $\l>\b$) of the problem (\ref{eq})-(\ref{c3}).  Such an assumption holds in the physical cases evoked in the introduction for the positive values of the parameter $m$.

When the proofs in the convex case are close to the ones of the concave case we will remove some details in order to shorten them.

\section{Concave solutions}
\begin{theorem}\label{concave}
Let $\a \in \mathbb R$ and $0\leq \l<\b$. If $g(x)<0$ for $x\in(\l,\b]$ and $g(\l)=0$, then the problem {\rm (\ref{eq})-(\ref{c3})} admits a unique concave solution.
\end{theorem}
\begin{proof}[Proof of existence] 
Let $f_\gamma$ be a solution of the initial value problem (\ref{initial}) with $\l<\b$ and $\gamma\leq 0$. As long as we have $f'_\gamma>\l$ and $f''_\gamma<0$, $f_\gamma$ exists. Because of Lemma \ref{vanish}, there are only three possibilities

\vspace{3mm}
(a) $f''_\gamma$ becomes positive from a point such that $f'_\gamma>\l,$

(b) $f'_\gamma$ takes the value $\l$ at some point for which $f''_\gamma<0,$

(c) we always have $f'_\gamma>\l$ and $f''_\gamma<0.$

\vspace{3mm}
\noindent As $f'_0(0)=\b>\l$, $f''_0(0)=0$ and $f'''_0(0)=-g(f'_0(0))=-g(\b)>0$ we have that $f'_0(t)>\l$ and $f''_0(t)>0$ on some interval $[0,t_0)$. By continuity it follows that $f''_\gamma$ becomes positive at some point for which $f'_\gamma>\l$ for small values of $-\gamma$ and thus $f_\gamma$ is of type (a).

On the other hand, as long as $f''_\gamma(t)<0$ and $f'_\gamma(t)\geq \l$, we have $f_\gamma(t)\geq \a$, $f'_\gamma(t)\leq \b$ and, using (\ref{i1}) we obtain
\begin{align*}
f''_\gamma(t) &\leq \g +|\a|\b+\a\b+\b^2t-\int_0^t g(f_\g'(\xi))d\xi \cr
&\leq \g +|\a|\b+\a\b+(\b^2+C)t
\end{align*}
where $C=\max\{-g(x)~;~x\in[\l,\b]\}>0$. Integrating once again we have
$$\l \leq f'_\gamma(t) \leq \frac{\b^2+C}{2}t^2+(\g+\a\b+|\a|\b)t+\b:=P_\g(t).$$
Hence, for $-\gamma$ large enough, the equation $P_\g(t)=\l$ has two positive roots $t_0<t_1$, and therefore we have $f'_\g(t_0)=\l$ and $f''_\gamma(t)<0$ for $t\leq t_0$, and $f_\gamma$ is of type (b).

Defining $A=\left \{ \gamma<0\, ; \, \text{$f_\gamma$ is of type (a)}\right \}$ and
$B=\left \{ \gamma<0\, ; \, \text{$f_\gamma$ is of type (b)}\right \}$ we have that 
$A \neq \emptyset$, $B \neq \emptyset$ and $A \cap B = \emptyset$. Both $A$ and $B$ are open sets, so there exists a $\gamma_*<0$ such that the solution $f_{\gamma_*}$ of (\ref{initial}) is of type (c) and is defined on the whole interval $[0,\infty)$. For this solution we have that $f'_{\gamma_*}>\l$ and $f''_{\gamma_*}<0$ which implies that $f'_{\gamma_*} \rightarrow l \in[\l,\b)$ as $t\rightarrow \infty$. From Lemma \ref{lim} and the fact that $g<0$ on $(\l,\b]$ we get $l=\l$.
\end{proof}

\vspace{3mm}
\begin{proof}[Proof of uniqueness]
Let $f$ be a concave solution of (\ref{eq})-(\ref{c3}). As $f'$ is positive and strictly decreasing, we can define a function $v:(\l^2,\b^2]\rightarrow[\a,\infty)$ such that
$$\forall t\geq 0,\quad v(f'(t)^2)=f(t).$$
Setting $y=f'(t)^2$ leads to
\begin{equation}
f(t)=v(y),\quad f''(t)=\frac{1}{2v'(y)}\quad \text{and} \quad f'''(t)=-\frac{v''(y)\sqrt{y}}{2v'(y)^3}.\label{y1}
\end{equation}
Then, using (\ref{eq}) we obtain
\begin{equation}
\forall y\in (\l^2,\b^2],\quad v''(y)=\frac{v(y)v'(y)^2}{\sqrt{y}}+\frac{2v'(y)^3g(\sqrt{y})}{\sqrt{y}} \label{eq-v1}
\end{equation}
with
$$v(\b^2)=v(f'(0)^2)=\a \quad  \text{and}\quad v'(\b^2)=\frac{1}{2\g}<0.$$
Suppose now that there are two concave solutions $f_{1}$ and $f_{2}$ of 
(\ref{eq})-(\ref{c3}) with $f''_i(0)=\gamma_i<0$ $i\in\{1,2\}$ and $\gamma_1>\gamma_2$. They gives $v_{1},\, v_{2}$ solutions of equation (\ref{eq-v1}) defined on
$(\l^2,\b^2]$ such that, if $w=v_1-v_2$, we have
$$w(\b^2)=0 \quad \text{and} \quad w'(\b^2)=\frac{1}{2\gamma_1}-\frac{1}{2\gamma_2}<0.$$
If $w'$ vanishes, there exists an $x$ in $(\l^2,\b^2]$ such that $w'(x)=0$, $w''(x)\leq0$ and $w(x)>0$. But from (\ref{eq-v1}) we then obtain
$$w''(x)=\frac{v'_1(x)^2}{\sqrt{x}}w(x)>0$$
and this is a contradiction. Therefore, $w'<0$ and $w>0$ on $(\l^2,\b^2]$. Set now $V_i=\frac{1}{v'_i}$ for $i\in\{1,2\}$ and $W=V_1-V_2$ to obtain
$$W'(y)=-\frac{w(y)}{\sqrt{y}}-\frac{2w'(y)g(\sqrt{y})}{\sqrt{y}}<0.$$
But, using (\ref{y1}) we have $V_i(f'_i(t)^2)=2f_i''(t)$ and thanks to Lemma \ref{f2} we get $W(y) \rightarrow 0$ as $y\rightarrow \l^2$. Since $W$ is decreasing and $W(\b^2)=2(\g_1-\g_2)>0$ this is a contradiction.
\end{proof}

\begin{remark} If, in addition to the hypotheses of Theorem $\ref{concave}$, the function $g$ is assumed to be non increasing, then we can write a much simpler proof for the uniqueness result.
For that, let $f_1$ and $f_2$ be two concave solutions of {\rm (\ref{eq})-(\ref{c3})} and let $\gamma_i=f''_i(0)<0$, $i\in \left\{ 1,2 \right \}$ with $\gamma_1>\gamma_2$. Writing $u=f_1-f_2$, we have $u'(0)=0$, $u'(\infty)=0$ and $u''(0)>0.$
Hence $u'$ admits a positive local maximum at some $t_0>0$ such that $u'(t)>0$ for $t\in(0,t_0]$. As $u$ is increasing on $[0,t_0]$ and $u(0)=0$ we have $u(t_0)>0$. Then, from $(\ref{eq})$ and since $f'_i>\l$, $f''_i<0$ and $f''_1(t_0)=f''_2(t_0)$ we get
\begin{equation*}
u'''(t_0)=-f''_1(t_0)u(t_0)+g(f_2'(t_0))-g(f_1'(t_0))>0
\end{equation*}
because $f_1'(t_0)>f_2'(t_0)>\l$ and a contradiction with the fact that $u'''(t_0)\leq 0$.
\end{remark}

The following Proposition gives some informations about the behaviour at infinity of the concave solution of the problem (\ref{eq})-(\ref{c3}) obtained in Theorem \ref{concave}.

\begin{proposition} Let $\a\in \R$ and $0\leq \l<\b$. Let us assume that $g<0$ on $(\l,\b]$ and $g(\l)=0$, and let $f$ be the concave solution of $(\ref{eq})$-$(\ref{c3})$. Then, there exists a constant $\mu$ such that $\a<\mu<\sqrt{\a^2+2(\b-\l)}$ and
\begin{equation}
\lim_{t\to\infty}\{f(t)-(\lambda t+\mu)\}=0. \label{p1}
\end{equation}
Moreover, for all $t\geq 0$, one has $\l t +\a \leq f(t) \leq \l t +\mu$.
\end{proposition}
\begin{proof} 
Since $f$ is concave, then for all $t\geq 0$ we
have
$f'(t)\in(\lambda,\beta]$ and the function $t\mapsto f(t)-\lambda t$ is
increasing. Hence $f(t)-\lambda t\to\mu\in(\alpha,\infty]$ as $t\to\infty$. In addition, we have
$$\forall t\geq 0,\quad f'''(t)=-f(t)f''(t)-g(f'(t))\geq-f(t)f''(t)$$
and thus 
$$\forall t\geq 0,\quad {f'''(t)\over -f''(t)}\geq f(t)\geq f(t)-\l t.$$
If we assume that $\mu=\infty$, it follows that
$$\lim_{t\to\infty}{f'''(t)\over -f''(t)}=\infty.$$
Therefore, there exists $t_0\geq 0$ such that $f'''(t)\geq -f''(t)$ for $t\geq t_0$.
Then integrating twice and using Lemma \ref{f2} we get 
$$\forall t\geq t_0,\quad -f''(t)\geq-\lambda+f'(t)$$
and 
$$\forall t\geq t_0,\quad -f'(t)+f'(t_0)\geq-\lambda t+\lambda t_0+f(t)-f(t_0).$$
Since the left hand side is bounded, we get a contradiction. Therefore, $\mu<\infty$ and we have
$$\forall t> 0,~~~~~~\lambda t+\a<f(t)<\lambda t+\mu.$$
Finally, let us introduce the auxiliary nonnegative function
$$u(t)=f'(t)+{1\over 2}(f(t)-\lambda t)^2.$$
From (\ref{p1}), we see that $u$ is bounded. Moreover, we have
$$u''(t)=-g(f'(t))-\lambda tf''(t)+(f'(t)-\lambda)^2>0$$
and $u$ is convex. Therefore $u$ is decreasing and thus
$$\b+{\a^2\over 2}=u(0)>u(\infty)=\l+{\mu^2 \over 2}.$$
This completes the proof.
\end{proof}

\begin{remark}
If $\l=0$ the previous result means that the concave solution of {\rm (\ref{eq})-(\ref{c3})} is bounded.
\end{remark}

We have the following nonexistence result
\begin{theorem} \label{nexist} 
Let $\alpha\leq 0$ and $0\leq\lambda<\beta$. If $g$ is differentiable and if
\begin{equation}
\forall x\in[\lambda,\beta],\quad g(x)\geq x^2-\lambda x\quad \text{and} \quad
-\alpha+\max_{x\in[\lambda,\beta]}\{g(x)-x^2+\lambda x\}>0, \label{g1}
\end{equation}
then the problem $(\ref{eq})$-$(\ref{c3})$ does not admit concave solutions.
\end{theorem}
\begin{proof}
We follow an idea of \cite{wang}. Let $0\leq \l<\b$ and suppose that $f$ is a concave solution of (\ref{eq})-(\ref{c3}). As $f'$ is positive and strictly decreasing, we can define a negative function $v:(\l,\b]\rightarrow {\mathbb R}$ such that
$$\forall t\geq 0,\quad v(f'(t))=f''(t).$$
Setting $y=f'(t)$ we obtain
\begin{equation}
f''(t)=v(y), \quad f'''(t)=v(y)v'(y) \quad \text{and} \quad f^{(4)}(t)=v(y)^2v''(y)+v(y)v'(y)^2.
\label{vy}
\end{equation}
Derivating equation (\ref{eq}) leads to
$$f^{(4)}+f'f''+ff'''+g'(f')f''=0$$
and as
$$f=\frac{-f'''-g(f')}{f''}$$
we get
\begin{equation}
v''(y)=-\frac{y}{v(y)}-\left ( \frac{g(y)}{v(y)} \right )'. \label{eqv}
\end{equation}
Integrating (\ref{eqv}) on $[z,z_1]$ with $\l\leq z \leq z_1$ leads to
\begin{equation}
v'(z_1)-v'(z)=-\int_z^{z_1} \frac{y}{v(y)}dy-\frac{g(z_1)}{v(z_1)}+\frac{g(z)}{v(z)}. \label{iv1}
\end{equation}
Using equation (\ref{eq}) and (\ref{vy}), equality (\ref{iv1}) becomes 
$$-f(s_1)-v'(z)=-\int_z^{z_1} \frac{y}{v(y)}dy+\frac{g(z)}{v(z)}$$
with $s_1$ such that $f'(s_1)=z_1$.

Integrating on $[\l,x]$ with $\l\leq x \leq z_1$, and since $v(\l)=f''(\infty)=0$ by Lemma \ref{f2}, we get 
\begin{align*}
-f(s_1)(x-\l)-v(x)&=-\int_\l^x\left ( \int_z^{z_1} \frac{y}{v(y)}dy \right )dz+\int_\l^x\frac{g(z)}{v(z)}dz \cr
&=-\left [ z \int_z^{z_1} \frac{y}{v(y)}dy \right ]_\l^x+ \int_\l^x \frac{-z^2+g(z)}{v(z)}dz \cr
&=-x \int_x^{z_1} \frac{y}{v(y)}dy+\l \int_\l^{z_1} \frac{z}{v(z)}dz + \int_\l^x \frac{-z^2+g(z)}{v(z)}dz
\end{align*}
and taking $x={z_1}$ and $z_1\to \b$ we derive
\begin{equation}
-v(\b)= \int_\l^\b \frac{\l z-z^2+g(z)}{v(z)}dz+\a(\b-\l). \label{int1}
\end{equation}
Since $v(\beta)\leq 0$, the right hand side of $(\ref{int1})$ must be nonnegative,
but this cannot be the case if $(\ref{g1})$ holds.
\end{proof}


\begin{remark}Using the previous Theorem we can recover the following nonexistence result
\begin{itemize}
\item for free convection $($i.e. $g(x)=-\frac{2m}{m+1}x^2$ and $\l=0$$)$ there is no concave solutions for $-1<m\leq -\frac{1}{3}$ when $\a<0$, and $-1<m< -\frac{1}{3}$ when $\a=0$ $($see {\rm \cite{brighicr}}, {\rm \cite{brighi02}}, {\rm \cite{brighi01}} and {\rm \cite{brighisari}}$)$,
\end{itemize}
and obtain the new results
\begin{itemize}
\item for the Falkner-Skan case $($i.e. $g(x)=m(1-x)(1+x)$ and $\l=1$$)$ there is no concave solutions with $\a\leq 0$ and $\b>1$ for $m\leq -\frac{\b}{1+\b}$,
\item for mixed convection $($i.e. $g(x)=\frac{2m}{m+1}x(1-x)$ and $\l=1$$)$ there is no concave solutions with $\a\leq 0$ and $\b>1$ for $-1<m< -\frac{1}{3}$.
\end{itemize}
\end{remark}

\section{Convex solutions}
\begin{theorem} \label{convex}
Let $\a\in\mathbb R$ and $0\leq \b<\l$. If $g(x)>0$ for $x\in[\b,\l)$ and $g(\l)=0$, then the problem {\rm (\ref{eq})-(\ref{c3})} admits a unique convex solution. 
\end{theorem}
\begin{proof}[Proof of existence] Let $f_\gamma$ be a solution of the initial value problem (\ref{initial}) with $0\leq \b<\l$ and $\gamma\geq 0$. We notice that $f_\gamma$ exists as long as we have $f''_\g>0$ and $f'_\g<\l$. From Lemma \ref{vanish}, $f''_\gamma$ cannot vanish at a point where $f'_\gamma=\l$ and it follows that there are only three possibilities

\vspace{3mm}
(a) $f''_\gamma$ becomes negative from a point such that $f'_\gamma<\l,$

(b) $f'_\gamma$ takes the value $\l$ at some point for which $f''_\gamma>0,$

(c) we always have $\b\leq f'_\gamma<\l$ and $f''_\gamma>0.$
\vspace{3mm}

\noindent As $f'_0(0)=\b<\l$, $f''_0(0)=0$ and $f'''_0(0)=-g(\b)<0$, we have that $f_0$ is of type (a), and by continuity it must be so for $f_\gamma$ with $\gamma>0$ small enough.

On the other hand, as long as $f''_\gamma(t)>0$ and $f'_\gamma(t)\leq \l$, we have  $f_\gamma(t)\leq \l t+\a$, and (\ref{i1}) leads to
\begin{align*}
f''_\gamma(t) &\geq \g-f_\gamma(t)f'_\gamma(t)+\a\b-\int_0^tg(f'_\gamma(\xi))d\xi \cr
&\geq \g-(\l t+|\a|)\l+\a\b-\int_0^tg(f'_\gamma(\xi))d\xi \cr
&\geq \g-|\a|\l+\a\b-(\l^2+C)t 
\end{align*}
where $C=\max\{g(x)~;~x\in[\b,\l]\}>0$ and integrating once again we have
\begin{equation*}
\l\geq f'_\gamma(t) \geq -\frac{\l^2+C}{2}t^2+(\g-|\a|\l+\a\b)t+\b:=P_\g(t).\label{p}
\end{equation*}
Hence, for $\gamma$ large enough, the equation $P_\g(t)=\l$ has two positive roots $t_0<t_1$, and therefore, for such a $\g$, we have $f'_\gamma(t_0)=\l$ and $f''_\gamma(t)>0$ for $t\leq t_0$, and $f_\gamma$ is of type (b).

Defining $A=\left \{ \gamma>0\, ; \, \text{$f_\gamma$ is of type (a)}\right \}$ and
$B=\left \{ \gamma>0\, ; \, \text{$f_\gamma$ is of type (b)}\right \}$ we have that 
$A \neq \emptyset$, $B \neq \emptyset$ and $A \cap B = \emptyset$. Both $A$ and $B$ are open sets, so there exists a $\gamma_*>0$ such that the solution $f_{\gamma_*}$ of 
(\ref{initial}) is of type (c) and is defined on the whole interval $[0,\infty)$. For this solution we have that $0<f'_{\gamma_*}<\l$ and $f''_{\gamma_*}>0$ which implies that $f'_{\gamma_*} \rightarrow l \in (\b,\l]$ as $t\rightarrow \infty$. From Lemma \ref{lim} and the fact that $g>0$ on $[\b,\l)$ we get $l=\l$.
\end{proof}

\vspace{3mm}
\begin{proof}[Proof of uniqueness]
Let $f$ be a convex solution of (\ref{eq})-(\ref{c3}). As $f'$ and $f''$ are positive, we can define a function $v:[\b^2,\l^2)\rightarrow[\a,\infty)$ such that
$$\forall t\geq 0,\quad v(f'(t)^2)=f(t).$$
We have
\begin{equation}
\forall y\in [\b^2,\l^2),\quad v''(y)=\frac{v(y)v'(y)^2}{\sqrt{y}}+\frac{2v'(y)^3g(\sqrt{y})}{\sqrt{y}} \label{eq-v}
\end{equation}
and
$$v(\b^2)=v(f'(0)^2)=\a \quad \text{and}\quad v'(\b^2)=\frac{1}{2\g}>0.$$
Suppose now that there are two convex solutions $f_{1}$ and $f_{2}$ of 
(\ref{eq})-(\ref{c3}) with $f''_i(0)=\gamma_i>0$, $i\in\{1,2\}$ and $\gamma_1>\gamma_2$. They gives $v_{1},\, v_{2}$ solutions of equation (\ref{eq-v}) defined on
$[\b^2,\l^2)$ such that
for $w=v_{1}-v_{2}$, we have $w(\b^2)=0$ and $w'(\b^2)< 0$. 
If $w'$ vanishes, there exists an $x$ in $[\b^2,\l^2)$ such that $w'(x)=0$, $w''(x)\geq0$ and $w(x)<0$. But from (\ref{eq-v}) we then obtain $w''(x)<0$ and a contradiction. Therefore $w'<0$ and $w<0$ on $[\b^2,\l^2)$. Set now $V_i=\frac{1}{v'_i}$ for $i\in\{1,2\}$ and $W=V_1-V_2$. We have $W>0$ and using (\ref{eq-v}) we obtain $W'(y)>0.$
But, $W(\b^2)=2(\g_1-\g_2)>0$ and $W(y) \rightarrow 0$ as $y\rightarrow \l^2$, this contradict the fact that $W$ is increasing.
\end{proof}

\begin{proposition} Let $\a\in \R$ and $0\leq \b<\l$. Assume that $g>0$ on $[\b,\l)$ and $g(\l)=0$, and let $f$ be the convex solution of $(\ref{eq})$-$(\ref{c3})$. Then, there exists a constant $\mu>\a$ such that 
$$\lim_{t\to\infty}\{f(t)-(\lambda t+\mu)\}=0.$$
Moreover, for all $t\geq 0$, one has $\l t +\a \leq f(t) \leq \l t +\mu$.
\end{proposition}
\begin{proof} 
Since $f$ is convex, for all $t\geq 0$ we have
$f'(t)\in[\beta,\lambda)$. Then the function $t\mapsto f(t)-\lambda t$ is decreasing
and thus $f(t)-\lambda t\to\mu\in[-\infty,\alpha)$ as $t\to\infty$. On the other
hand, we have
$$\forall t\geq 0, \quad f'''(t)=-f(t)f''(t)-g(f'(t))\leq-f(t)f''(t)$$
and since $f(t)\to\infty$ as $t\to\infty$,
there exists $t_0\geq 0$ such that $f'''(t)\leq -f''(t)$ for $t\geq t_0$.
Then integrating twice and using Lemma \ref{f2} we get 
$$\forall t\geq t_0,\quad -f''(t)\leq-\lambda+f'(t)$$
and 
$$\forall t\geq t_0,\quad -f'(t)+f'(t_0)\leq-\lambda t+\lambda t_0+f(t)-f(t_0).$$
Since the left hand side is bounded, we necessarily get $\mu>-\infty$, and then we have $\l t +\mu \leq f(t) \leq \l t +\a$ for all $t\geq 0$.
\end{proof}

\bigskip
Let us finish this section with the following nonexistence result.
\begin{theorem}
Let $\alpha\leq 0$ and $0\leq\b<\l$. If $g$ is differentiable and if
\begin{equation}
\forall x\in[\b,\l],\quad g(x)\leq x^2-\lambda x\quad \text{and} \quad
-\alpha+\max_{x\in[\b,\l]}\{x^2-\lambda x-g(x)\}>0, \label{g2}
\end{equation}
then the problem $(\ref{eq})$-$(\ref{c3})$ does not admit convex solutions
\end{theorem}
\begin{proof}
Let $0\leq \b<\l$ and suppose that $f$ is a convex solution of (\ref{eq})-(\ref{c3}). As $f'$ and $f''$ are positive, we can define a positive function $v:[\b,\l)\rightarrow {\mathbb R}$ such that
$$\forall t\geq 0,\quad v(f'(t))=f''(t).$$
Following the method used in the proof of Theorem \ref{nexist}, we see that this function $v$ satisfies (\ref{eqv}) and 
\begin{equation}
-v(\b)= \int_\b^\l \frac{y^2-\l y-g(y)}{v(y)}dy-\a(\l-\b). \label{eqa}
\end{equation}
Since $v(\beta)\geq 0$, the right hand side of $(\ref{eqa})$ must be nonpositive,
but this cannot be the case if $(\ref{g2})$ holds.
\end{proof}
\bigskip

\begin{remark}Using the previous Theorem we can recover the following nonexistence result
\begin{itemize}
\item for the Falkner-Skan case $($i.e. $g(x)=m(1-x)(1+x)$ and $\l=1$$)$ there is no convex solutions for $m\leq -\frac{1}{2}$ $($see {\rm \cite{wang}}$)$,
\end{itemize}
and obtain the new result
\begin{itemize}
\item for mixed convection $($i.e. $g(x)=\frac{2m}{m+1}x(1-x)$ and $\l=1$$)$ there is no convex solutions for $-1<m\leq -\frac{1}{3}$.
\end{itemize}
\end{remark}


\section{Conclusion}
In this paper we have obtained existence, uniqueness and nonexistence results for the concave or convex solutions of a general boundary value problem arising in many fields of application under some reasonable hypotheses. 
All these hypotheses are verified in important physical cases in the framework of boundary layer approximations such as free or mixed convection
, flow adjacent to stretching walls, high frequency excitation of liquid metal and two dimensional flow of a slightly viscous incompressible fluid past a wedge. 


All our results hold for $\b$, $\l$ and $g$ such that the function $g$ vanishes at $\l$ but does not vanish between $\b$ and $\l$. 
Of course, under the same hypotheses, 
solutions whose concavity changes may exist as it can be seen in \cite{brighisari} for the case of free convection and in \cite{hast1} for the Falkner-Skan problem when  the parameter $m$ is positive.

In addition, if we assume that the sign of $g$ is the opposite of the one that we have in the Theorems \ref{concave} and \ref{convex}, then we can have multiple concave or convex solutions, as it is the case for free or mixed convection and for the Falkner-Skan problem when  the parameter $m$ is negative. See for example \cite{brighi01}, \cite{gaeta}, \cite{guedda}, \cite{guedda2} and \cite{hart}.



\end{document}